\documentclass[11pt]{amsart}%{article}

\usepackage{amsfonts}
\usepackage{mathtools}
\usepackage{blindtext}
\usepackage{tikz}
\usepackage{amsmath}
\usepackage{graphicx}

\newtheorem{thm}{Theorem}

\newtheorem{prop}{Proposition}

\newtheorem{remark}{Remark}

\PassOptionsToPackage{usenames,dvipsnames}{xcolor}
\mathtoolsset{showonlyrefs}

\newdimen\slantmathcorr
\def\oversl#1{\setbox0=\hbox{$#1$}
\slantmathcorr=\wd0
\hskip 0.2\slantmathcorr \overline{\hbox to 0.8\wd0{\vphantom{\hbox{$#1$}}}}
\hskip-\wd0\hbox{$#1$}
}
\def\undersl#1{\setbox0=\hbox{$#1$}
\slantmathcorr=\wd0
\underline{\hbox to 0.8\wd0{\vphantom{\hbox{$#1$}}}}
\hskip-0.8\wd0\hbox{$#1$}
}

\newcommand{\be}{\begin{equation}}
\newcommand{\ee}{\end{equation}}
\newcommand{\ben}{\begin{equation*}}
\newcommand{\een}{\end{equation*}}
\newcommand{\ba}{\begin{aligned}}
\newcommand{\ea}{\end{aligned}}

\newcommand{\gep}{\epsilon}

\newcommand{\whH}{\widehat{H}}

\newcommand{\whL}{\widehat{L}}

\newcommand{\whk}{\widehat{\kappa}}

\newcommand{\whq}{\widehat{q}}

\newcommand{\pibar}{\overline{\Pi}}

\newcommand{\ga}{\alpha}
\newcommand{\gl}{\lambda}
\newcommand{\gk}{\kappa}

\newcommand{\gd}{\delta}

\newcommand{\vep}{\varepsilon}

\newcommand{\whn}{\widehat{n}}

\newcommand{\tx}{\tau_x}

\numberwithin{equation}{section}
\numberwithin{thm}{section}
\numberwithin{lemma}{section}
\numberwithin{prop}{section}
\numberwithin{remark}{section}
\numberwithin{cor}{section}
\numberwithin{ex}{section}

\begin{document}

\title{Cram\'er's estimate for the reflected process revisited}

%    author one information
\author{R. A. Doney}
\address{R. A. Doney}
\curraddr{School of Mathematics, University of Manchester, Oxford Road, Manchester,  M13 9PL, UK}
\email{ron.doney@manchester.ac.uk}
%\thanks{This is the second author}

%    author two information
\author{Philip S. Griffin}
\address{Philip S. Griffin}
\curraddr{Department of Mathematics, Syracuse University, Syracuse, NY 13244-1150, USA}
\email{psgriffi@syr.edu}
\thanks{This work was partially supported by Simons Foundation Grant
226863 to Philip Griffin}


\maketitle
\begin{abstract}
The reflected process of a random walk or L\'{e}vy process arises in many
areas of applied probability, and a question of particular interest is how
the tail of the distribution of the heights of the excursions away from zero
behaves asymptotically. The L\'{e}vy analogue of this is the tail behaviour of the characteristic measure of the height of an excursion.
Apparently the only case where this is known is when Cram\'{e}r's condition
hold. Here we establish the asymptotic behaviour for a large class of L\'{e}%
vy processes which have exponential moments but do not satisfy Cram\'{e}r's
condition. Our proof also applies in the Cram\'{e}r case, and corrects a
proof of this given in Doney and Maller \cite{DM}.
\end{abstract}

%\date{}

%%%%%%%%%%%%%%%%%%%%%%%%%%%%%%%%%%%%%%%%%%%%%%%%%%%%%%%%%%%%%%%%%%%%%%%%%%%%%%%%%%%%%%%%%%%%%%%%%%%%%%%%%%%
%Section Introduction
%%%%%%%%%%%%%%%%%%%%%%%%%%%%%%%%%%%%%%%%%%%%%%%%%%%%%%%%%%%%%%%%%%%%%%%%%%%%%%%%%%%%%%%%%%%%%%%%%%%%%%%%%%%

\setcounter{equation}{0}

\section{Introduction}

\label{Intro}

The reflected process $R=(R_{n},n\geq 0)$ formed from a random walk $%
S=(S_{n},n\geq 0)$ by setting%
\begin{equation*}
R_{n}=S_{n}-I_{n}\text{ where }I_{n}=\min_{i\leq n}S_{i},\text{ }n\geq 0,
\end{equation*}%
arises in many areas of applied probability, including queuing theory, risk
theory, and mathematical genetics. In all these areas, the i.i.d sequence of
random variables defined by%
\begin{equation*}
h_{i}=\max_{0\leq n\leq \hat{T}_{i}-\hat{T}_{i-1}}\{S_{\hat{T}_{i-1}+n}-S_{%
\hat{T}_{i-1}}\},\text{ }i=1,2,\cdots ,
\end{equation*}%
where $\hat{T}_{i}$ is the $i$th strict descending ladder time, with $\hat{T}%
_{0}=0,$ is of central importance. %interest. 
These random variables give the heights of the excursions of $R$ away from $%
0 $, or equivalently the heights of the excursions of $S$ above its minimum.
Our main focus will be on the asymptotic behaviour of $P(h_{1}>x)$ which
among other things is useful in the study of the point process of excursion
heights.

In continuous time we replace the random walk by a L\'{e}vy process $%
X=(X_{t},t\geq 0)$ and study 
\begin{equation}  \label{R}
R=(R_{t},t\geq 0)\text{ with }R_{t}=X_{t}-\setbox0=\hbox{$X$} \slantmathcorr=%
\wd0 \underline{\hbox to 0.8\wd0{\vphantom{\hbox{$X$}}}} \hskip-0.8\wd0%
\hbox{$X$} _{t},\text{ where }\setbox0=\hbox{$X$} \slantmathcorr=\wd0 
\underline{\hbox to 0.8\wd0{\vphantom{\hbox{$X$}}}} \hskip-0.8\wd0\hbox{$X$}
_{t}=\inf_{s\leq t}X_{s}.
\end{equation}%
In mathematical finance $R$ is called the drawup. When indexed by local time
at the infimum, the excursions of $R$ away from $0$ form a Poisson point
process whose characteristic measure we denote by $\hat{n}.$ If $h$ denotes
the height of a generic excursion, then $\hat{n}(h>x)$ is the L\'evy
analogue of $P(h_{1}>x)$.

Our main interest is in the L\'{e}vy process case, but we start by reviewing
some discrete time results. A classical case where the asymptotic behaviour
of $P(h_{1}>x)$ is known is when $S$ satisfies Cram\'{e}r's condition,
namely $E(e^{\gamma S_{1}})=1$ for some $\gamma \in (0,\infty ).$ Then $S$
drifts to $-\infty $ and for $x>0$ the first time passage $\tau
_{x}=\inf\{n:S_n>x\}$ of $S$ to $(x,\infty )$ is defective and satisfies%
\begin{equation}  \label{a}
\lim_{x\rightarrow \infty }e^{\gamma x}P(\tau _{x}<\infty )=\Gamma ,
\end{equation}%
where $\Gamma $ is a known finite constant. It then follows immediately from
the identity%
\begin{equation}  \label{b}
P(\tau _{x}<\infty )=P(h_{1}>x)+\int_{0}^{\infty }P(h_{1}\leq x,|S_{\hat{T}%
_{1}}|\in dy)P(\tau _{x+y}<\infty ),
\end{equation}%
where $|S_{\hat{T}_{1}}|$ is the first strict descending ladder height, that%
\begin{equation}  \label{c}
\lim_{x\rightarrow \infty }e^{\gamma x}P(h_{1}>x)=\Gamma \{1-E(e^{-\gamma
|S_{\hat{T}_{1}}|})\}.
\end{equation}%
This argument is due to Iglehart \cite{Ig}.

The L\'{e}vy version of Cram\'{e}r's condition is that $E(e^{\gamma
X_{1}})=1 $ for some $\gamma >0$, and assuming this Bertoin and Doney \cite%
{BD94} proved the following analogue of \eqref{a}; 
\begin{equation}
\lim_{x\rightarrow \infty }e^{\gamma x}P(\tau _{x}<\infty )=\Gamma ^{\ast },
\label{h}
\end{equation}%
where $\tau _{x}=\inf\{t:X_t>x\}$ is now the first time passage of $X$ to $%
(x,\infty )$ and $\Gamma ^{\ast }$ is a known finite constant. The analogue
of \eqref{c} now becomes%
\begin{equation}  \label{g}
\lim_{x\rightarrow \infty }e^{\gamma x}\hat{n}(h>x)=\hat{\kappa}(\gamma
)\Gamma ^{\ast },
\end{equation}%
where $\hat{\kappa}$ is the Laplace exponent of the strictly decreasing
ladder height process. A proof of this result was given in \cite{DM}, but
there is a problem with the argument presented there. Specifically, equation
(15) therein is not fully justified, and we believe that it cannot be
justified. So our first aim is to rectify this, and we do so by using a
different approach which applies to a much more general situation.

For any non-negative function $f,$ let us say that $f\in \mathcal{L}%
^{(\alpha)},\alpha \geq 0,$ if%
\begin{equation}  \label{d}
\lim_{x\rightarrow \infty }\frac{f(x+y)}{f(x)}=e^{-\alpha y}\text{ for all }%
y.
\end{equation}%
For a random variable, $Z\in \mathcal{L}^{(\alpha)}$ means $P(Z>x)\in 
\mathcal{L}^{(\alpha) },$ and for a measure $\mu \in \mathcal{L}^{(\alpha)}$
means that $\overline{\mu }(x):=\mu ((x,\infty ))\in \mathcal{L}^{(\alpha)}.$
So for a random walk in the Cram\'{e}r case, if $\Gamma\not =0$, then $%
P(\tau _{x}<\infty )\in \mathcal{L}^{(\gamma) },h_{1}\in \mathcal{L}%
^{(\gamma )},$ and the ratio of $P(h_1>x)$ to $P(\tau_x<\infty)$ converges
to the constant $1-E(e^{-\gamma \hat{H}_{1}})$. Our first main result
includes the L\'{e}vy process version of this fact, but much more as well.

\begin{thm}
\label{T1} Fix $\alpha >0$. For any L\'evy process $X$, 
\begin{equation}  \label{tr}
P(\tau_x<\infty)\in \mathcal{L}^{(\alpha)}
\end{equation}
if and only if 
\begin{equation}  \label{er}
\widehat{n}(h>x)\in\mathcal{L}^{(\alpha)}
\end{equation}
in which case 
\begin{equation}  \label{ter}
\lim_{x\to\infty} \frac{\widehat{n}(h>x)}{P(\tau_x<\infty)}= \widehat{\kappa}%
(\alpha).
\end{equation}
\end{thm}

Thus in particular, \eqref{g} is now proved provided $\Gamma^*\not = 0$. (If 
$\Gamma^* = 0$, \eqref{g} continues to hold; see Remark \ref{Rem0}). Since
distributions in $\mathcal{L}^{(\alpha)}$ are ``close to exponential", this
result will also lead to useful Cram\'{e}r-like estimates for $\hat{n}(h>x)$
if we can replace the condition \eqref{tr} by a condition expressed in terms
of $\Pi_X ,$ the L\'{e}vy measure of $X$. We will give a complete answer to
this under the natural assumption that $\overline{\Pi}_X\in \mathcal{L}%
^{(\alpha)},$ but first we consider the situation that $\overline{\Pi}_X\in 
\mathcal{S}^{(\alpha)}$, the class of $a-$convolution equivalent functions,
for some $\alpha >0$. This means that $\overline{\Pi }_X\in \mathcal{L}%
^{(\alpha)},$ and additionally the probability distribution defined by $%
G(dy)=\Pi_X (dy)/\overline{\Pi }_X(1)$ for $y\in (1,\infty )$ satisfies%
\begin{equation}
\lim_{x\rightarrow \infty }\frac{\overline{G\ast G}(x)}{2\overline{G}(x)}%
=\int_{1}^{\infty }e^{\alpha y}G(dy)<\infty .  \label{m}
\end{equation}%
In this scenario $E(e^{\alpha X_{1}})<\infty,$ and since \eqref{tr} implies $%
X_t\to -\infty$ a.s., we can then assume, WLOG, that $E(e^{\alpha X_{1}})<1.$
This is because if not there exists a $\gamma\in(0,\alpha]$ such that $%
E(e^{\gamma X_{1}})=1,$ so we are back in the Cram\'{e}r situation. When $%
\overline{\Pi }_X\in \mathcal{S}^{(\alpha)}$ and $E(e^{\alpha X_{1}})<1,$ it
has been shown in Kl\"{u}ppelberg, Kyprianou and Maller \cite{KKM}, Lemma
3.5, that%
\begin{equation}  \label{n}
\lim_{x\rightarrow \infty }\frac{P(\tau _{x}^{H}<\infty )}{\overline{\Pi }%
_{H}(x)}=\frac{q}{\kappa (-\alpha )^{2}},
\end{equation}%
where $\Pi _{H}$ is the L\'{e}vy measure and $\tau _{x}^{H}$ the first
passage time for the increasing ladder height process $H$, and $\kappa $ and 
$q$ are the Laplace exponent and killing rate of $H$ respectively. Since $%
P(\tau _{x}^{H}<\infty )=P(\tau _{x}<\infty )$ and it was also claimed in
Proposition 5.3 of \cite{KKM} that $\overline{\Pi }_X\in \mathcal{L}%
^{(\alpha)}$ if and only if $\overline{\Pi }_{H}\in \mathcal{L}^{(\alpha)}$
and then $\overline{\Pi }_X(x)\sim \hat{\kappa}(\alpha )\overline{\Pi }%
_{H}(x),$ \eqref{n} is apparently equivalent to%
\begin{equation}  \label{o}
\lim_{x\rightarrow \infty }\frac{P(\tau _{x}<\infty )}{\overline{\Pi }_X(x)}=%
\frac{q}{\widehat{\kappa}(\alpha )\kappa (-\alpha )^{2}}.
\end{equation}%
Together with our Theorem \ref{T1} this would solve the problem in this
convolution equivalent case. However there is a problem with the proof of the
claimed equivalence of $\overline{\Pi }_X$ and $\overline{\Pi }_{H}$, specifically in
display (7.18) of \cite{KKM}, which we circumvent in proving

\begin{thm}
\label{KKMC} Fix $\alpha >0$. For any L\'evy process $X$, 
\begin{equation}  \label{piXr}
\overline{\Pi}_X\in \mathcal{L}^{(\alpha)}
\end{equation}
if and only if 
\begin{equation}  \label{piHr}
\overline{\Pi}_H\in \mathcal{L}^{(\alpha)}
\end{equation}
in which case 
\begin{equation}  \label{piXHr}
\lim_{x\to\infty} \frac{\overline{\Pi}_X(x)}{\overline{\Pi}_{H}(x)}= 
\widehat{\kappa}(\alpha).
\end{equation}
\end{thm}

\begin{remark}
Note that, unlike Prop 5.3 of \cite{KKM}, we do not require the assumption
that $X_{t}\rightarrow -\infty $ a.s. in this result.
\end{remark}

Our last main result addresses the possibility that there are situations
where $\overline{\Pi }_X\in \mathcal{L}^{(\alpha) }\backslash \mathcal{S}%
^{(\alpha) }$ and $P(\tau _{x}<\infty )$ (and hence $\hat{n}(h>x))$ has the
same asymptotic behaviour as $\overline{\Pi }_X(x).$

\begin{thm}
\label{LPP} Assume $\alpha >0,$ $\overline{\Pi }_X\in \mathcal{L}^{(\alpha)}$
and $E(e^{\alpha X_{1}})<1.$ Then 
\begin{equation}  \label{limXa0}
\lim_{x\rightarrow \infty }\frac{P(\tau _{x}<\infty )}{\overline{\Pi }_X(x)}%
=L\in (0,\infty )
\end{equation}%
if and only if $\overline{\Pi }_X\in \mathcal{S}^{(\alpha)}.$ In this case $%
\displaystyle L=\frac{q}{\widehat{\kappa}(\alpha )\kappa (-\alpha )^{2}},$
and%
\begin{equation}  \label{ter1}
\lim_{x\to\infty} \frac{P(\tau_x<\infty)}{\overline{\Pi}_{H}(x)}=\frac{q}{%
\kappa(-\alpha)^2}.
\end{equation}
\end{thm}

\begin{remark}
\label{Rem4} Note that the assumptions are equivalent to $\overline{\Pi }%
_{H}\in \mathcal{L}^{(\alpha)}$ , $E(e^{\alpha H_{1}})<1,$ and%
\begin{equation}  \label{tHrat}
\lim_{x\rightarrow \infty }\frac{P(\tau _{x}^{H}<\infty )}{\overline{\Pi }%
_{H}(x)}=L^{\prime }\in (0,\infty )
\end{equation}%
and because of Theorem \ref{KKMC}, the conclusion can be written as $%
\overline{\Pi }_{H}\in \mathcal{S}^{(\alpha)}$ and $L^{\prime }=q/\kappa
(-\alpha )^{2}. $ In fact our proof shows that this version of the result
holds for any defective subordinator.
\end{remark}

\begin{remark}
Note that, in particular, our results show that when $\alpha >0,$ $\overline{%
\Pi }_{X}\in \mathcal{S}^{(\alpha )}$ and $E(e^{\alpha X_{1}})<1$ the
quantities $\overline{\Pi }_{H}(x),P(\tau _{x}<\infty )$ and $\widehat{n}%
(h>x)$ all have the same asymptotic behaviour as $\overline{\Pi }_{X}(x).$
This contrasts with the Cram\'{e}r case, when $P(\tau _{x}<\infty )$ and $%
\widehat{n}(h>x)$ are comparable to each other but not to $\overline{\Pi }_{X}(x)$  since then $%
\overline{\Pi }_{X}(x)=o(e^{-\gamma x}).$
\end{remark}

We conclude this section by remarking that exactly analogous results hold in
the discrete time setting, and their proofs, which we omit, are considerably
simpler. Also, our techniques yield some results for the case $\alpha =0.$
These can be found in the remarks in Section \ref{sPro}.

%%%%%%%%%%%%%%%%%%%%%%%%%%%%%%%%%%%%%%%%%%%%%%%%%%%%%%%%%%%%%%%%%%%%%%%%%%%%%%%%%%%%%%%%%%%%%%%%%%%%%%%%%%%
%Section Preliminaries
%%%%%%%%%%%%%%%%%%%%%%%%%%%%%%%%%%%%%%%%%%%%%%%%%%%%%%%%%%%%%%%%%%%%%%%%%%%%%%%%%%%%%%%%%%%%%%%%%%%%%%%%%%%

\setcounter{equation}{0}

\section{ Preliminaries}

\label{sPre}

We briefly collect the pertinent properties of a L\'evy process to be used
in this paper. Further details can be found for example in \cite{bert}, \cite%
{D}, \cite{kypbook} and \cite{sato}. Let $(L^{-1}_s,H_s)_{s \geq 0}$ denote
the \textit{weakly} 
%\footnote{The distinction between weak and strict only makes a difference when $X$ is compound Poisson.}
ascending bivariate ladder process of $X$. When $X_t\to -\infty$ a.s., $%
(L^{-1},H)$ is defective and may be obtained from a nondefective process by
exponential killing at some appropriate rate $q > 0$. When the process is
killed it is sent to some cemetery state, in which case probabilities and
expectation are understood to be taken over only non cemetery values. The
renewal function of $H$ is 
\begin{equation}  \label{Vkdef}
V(x)= \int_0^\infty P(H_s\le x)ds.
\end{equation} 
Note that 
$V(\infty):=\lim_{x\to\infty} V(x)=q^{-1}. $ 
The Laplace exponent $\kappa$ of $H$, defined by 
$e^{-\kappa(\lambda)} = E e^{ -\lambda{H}_1} $ 
for values of $\lambda\in \mathbb{R}$ for which the expectation is finite,
satisfies 
\begin{equation}
\kappa(\lambda)=q+d\lambda+\int_0^\infty (1-e^{-\lambda x})\Pi_H(dx).
\end{equation}
Observe that 
\begin{equation}  \label{Vk}
\int_{y\ge 0} e^{-\lambda y}V(dy)=\frac{1}{\kappa(\lambda)}
\end{equation}
for all $\lambda\in \mathbb{R}$ with $\kappa(\lambda)>0$.

Let $\widehat{X}_t=-X_t$, $t\ge 0$, denote the dual process, and $(\widehat{L}%
^{-1}, \widehat{H})$ the corresponding \textit{strictly} ascending bivariate
ladder processes of $\widehat{X}$. All quantities relating to $\widehat{X}$
will be denoted in the obvious way, for example $\widehat{\kappa}, \widehat{d%
}, \Pi_{\widehat{H}}$ and $\widehat{V}$. We may assume the normalisations of 
$L$ and $\widehat{L}$ are chosen so that the constant in the Wiener-Hopf
factorisation is 1; see (4) in Section VI.2 of \cite{bert}. 
$\widehat{L}$ is a local time at $0$ for the reflected process $R$, % 
and the excursion $e_t$ of $R$ at local time $t$ is given by 
\begin{equation*}
e_t(s)=X_{(\widehat{L}_{t-}^{-1}+s)\wedge\widehat{L}_{t}^{-1}}-\setbox0=%
\hbox{$X$} \slantmathcorr=\wd0 \underline{\hbox to
0.8\wd0{\vphantom{\hbox{$X$}}}} \hskip-0.8\wd0\hbox{$X$} _{\widehat{L}%
_{t-}^{-1}}.
\end{equation*}
If $e_t\not\equiv 0$, that is $\Delta \widehat{L}_{t}^{-1}>0$, then $e_t$
takes values in the space of excursions 
\begin{equation}
\mathcal{E}=\{\epsilon\in D:\epsilon(s)\ge 0 \text{ for all } 0\le s<\zeta,\
\zeta>0\},
\end{equation}
where $D$ is the Skorohod space of cadlag functions and $\zeta=\zeta(%
\epsilon)=\inf\{s:\epsilon(u)=\epsilon(v) \text{ all } u,v\ge s\}$ is the
lifetime of the excursion. Futhermore, $\{(t,e_t):e_t\in \mathcal{E}\}$ is a
Poisson point process with intensity (excursion) measure $\widehat{n}$.

For $\epsilon\in \mathcal{E}$, let $h=h(\epsilon)=\sup_{s\ge 0} \epsilon(s)$
be the height of the excursion $\epsilon$. Note that $\widehat{n}(h=0)>0$ if
and only if $X$ is compound Poisson. Set $|\widehat{n}|=\widehat{n}(\mathcal{%
E})=\widehat{n}(h\ge 0)$. The following result describing when $\widehat{n}$
is finite will be needed in the proof of Theorem \ref{T1}.

\begin{prop}
\label{P11} The excursion measure $\widehat{n}$ is finite if and only if one
of the following two conditions hold 
\begin{equation}  \label{nfi1}
0 \text{ is irregular for $[0,\infty)$ and $\overline{\Pi}_X^+(0)<\infty$};
\end{equation}
\begin{equation}  \label{nfi2}
0 \text{ is irregular for $(-\infty,0)$.}
\end{equation}
\end{prop}

\noindent\textbf{Proof of Proposition \ref{P11}}\ \ Excursion intervals are
precisely the non-empty intervals of the form $(\widehat{L}^{-1}_{t-}, 
\widehat{L}^{-1}_{t})$. Let 
\begin{equation}
T=\inf\{t:\Delta \widehat{L}^{-1}_{t}>0\}.
\end{equation}
Then $|\widehat{n}|=\infty$ iff $T=0$ a.s. We consider the three possible
cases;\medskip

\noindent Case I: $0$ is regular for both $[0,\infty)$ and $(-\infty, 0)$:

Then there are excursion intervals with end points arbitrarily close to $0$,
i.e. there exist $t_n\downarrow $ such that $\Delta\widehat{L}^{-1}_{t_n}>0$
and $\widehat{L}^{-1}_{t_n}\to 0$. If $t_n\downarrow s>0$ then by right
continuity, $\widehat{L}^{-1}_{s}= 0$. This implies $\widehat{L}^{-1}$ is
compound Poisson which is impossible when $0$ is regular for $(-\infty, 0)$.
Thus $T=0$ and so $|\widehat{n}|=\infty$.

In the two remaining cases, $0$ is irregular for exactly one of $[0,\infty)$
or $(-\infty, 0)$. In particular this implies $X$ has bounded variation and
so $X_t=Y_t-Z_t+ct$ where $X$ and $Y$ are pure jump subordinators.\medskip

\noindent Case II: $0$ is irregular for $[0,\infty)$:

In this case $c\le 0$ and $\widehat{L}^{-1}$ is not compound Poisson. Let 
\begin{equation}
S=\inf\{s:\Delta X_s>0\}.
\end{equation}
Then $S=\widehat{L}^{-1}_{T-}$ where $\widehat{L}^{-1}_{0-}=0$. By right
continuity of $\widehat{L}^{-1}$, $|\widehat{n}|<\infty$ precisely when $S>0$
a.s. which in turn is equivalent to $\overline{\Pi}_X^+(0)<\infty$.\medskip

\noindent Case III: $0$ is irregular for $(-\infty,0)$:

In this case $\widehat{L}^{-1} $ is compound Poisson by construction, see
p24 of \cite{D}, and so $T>0$. Thus $|\widehat{n}|<\infty$. \qquad$%
\sqcup\!\!\!\!\sqcap$ \bigskip

%%%%%%%%%%%%%%%%%%%%%%%%%%%%%%%%%%%%%%%%%%%%%%%%%%%%%%%%%%%%%%%%%%%%%%%%%%%%%%%%%%%%%%%%%%%%%%%%%%%%%%%%%%%
%Section Lalpha
%%%%%%%%%%%%%%%%%%%%%%%%%%%%%%%%%%%%%%%%%%%%%%%%%%%%%%%%%%%%%%%%%%%%%%%%%%%%%%%%%%%%%%%%%%%%%%%%%%%%%%%%%%%

\setcounter{equation}{0}

\section{$\mathcal{L}^{(\protect\alpha)}$ and $\mathcal{S}^{(\protect\alpha%
)} $}
\label{sLa}

Assume $f:(0,\infty)\to(0,\infty)$ satisfies 
\begin{equation}  \label{flim}
\lim_{x\to\infty}\frac{f(x+y)}{f(x)} \text{ exists for all } y>0.
\end{equation}
Then $g(x)=f(\ln x)$ is regularly varying at infinity with some index $%
-\alpha$ and hence 
\begin{equation}  \label{Ladef}
\lim_{x\to\infty}\frac{f(x+y)}{f(x)}=e^{-\alpha y} \text{ for all } y.
\end{equation}
Thus \eqref{Ladef} is equivalent to
the seemingly weaker \eqref{flim}. Exploiting the connection with regularly
varying functions further, 
a very useful global bound for the ratio in \eqref{Ladef} can be obtained
from Potter's Theorem. By applying Theorem 1.5.6(ii) of \cite{BGT} to the
function $g(x)=(x\vee e)^\ga f(\ln(x\vee e))$, it follows that if $f\in \mathcal{L%
}^{(\alpha)}$ is bounded away from $0$ and $\infty$ on compact subsets of $%
[1,\infty)$, then for every $\varepsilon>0$ there exists an $A=A_\vep$ such
that 
\begin{equation}  \label{Pot}
\frac{f(x+y)}{f(x)}\le A\left(e^{-(\alpha-\varepsilon)y}\vee
e^{-(\alpha+\varepsilon)y}\right)\ \text{ for all $x\ge 1, y\ge 1-x$}.
\end{equation}

The definition of $\overline{\Pi}_X\in\mathcal{S}^{(\alpha)}$ for $\alpha>0$
given in the introduction, applies equally well when $\alpha=0$. Here we
give an slightly different formulation which will be used later. Let $%
Z_1,Z_2 $ be independent and distributed as $Z$. Then $Z\in \mathcal{S}%
^{(\alpha)}, \alpha\ge 0$, if $Z\in\mathcal{L}^{(\alpha)}$ and 
\begin{equation}  \label{Sa}
\lim_{x\to\infty}\frac{P(Z_1+Z_2>x)}{P(Z_1>x)}\ \text{ exists.}
\end{equation}
Thus there is no requirement on the value of the limit in \eqref{Sa}.
However, see for example the discussion in Section 5 of \cite{W}, in this
case $Ee^{\alpha Z}<\infty$ and the limit in \eqref{Sa} is given by $%
2Ee^{\alpha Z}$. Thus $\overline{\Pi}_X\in \mathcal{S}^{(\alpha)}$, $%
\alpha\ge 0$, if $Z\in \mathcal{S}^{(\alpha)}$ where $Z$ has distribution
given by 
\begin{equation}
P(Z\in dy)=\frac{I(y>1)\Pi_X(dy)}{\overline{\Pi}_X(1)}.
\end{equation}
Since $\mathcal{S}^{(\alpha)}$ and $\mathcal{L}^{(\alpha)}$ are both closed
under tail equivalence the choice of cut-off point is not important.

%%%%%%%%%%%%%%%%%%%%%%%%%%%%%%%%%%%%%%%%%%%%%%%%%%%%%%%%%%%%%%%%%%%%%%%%%%%%%%%%%%%%%%%%%%%%%%%%%%%%%%%%%%%
%Section Proofs
%%%%%%%%%%%%%%%%%%%%%%%%%%%%%%%%%%%%%%%%%%%%%%%%%%%%%%%%%%%%%%%%%%%%%%%%%%%%%%%%%%%%%%%%%%%%%%%%%%%%%%%%%%%

\setcounter{equation}{0}

\section{Proofs}

\label{sPro}

Applying Corollary 4.1 of \cite{G} to the dual process $\widehat{X}$, the
L\'evy measure of $\widehat{H}$ is related to $\widehat{n}$ by the formula 
\begin{equation}  \label{piHn}
\Pi_{\widehat{H}}(dx)=\widehat{n}(|\epsilon(\zeta)|\in dx)+d_{{\widehat{L}}%
^{-1}}\Pi_X^-(dx),\ \ x>0,
\end{equation}
where $\Pi_X^-((x,\infty))=\Pi_X((-\infty,-x))$ for $x>0$. The final term on
the right hand side allows for the possibility of $X$ jumping down from a
strict current minimum. It is only present when $d_{{\widehat{L}}^{-1}}>0$,
which in turn implies $X$ has bounded variation. The Poisson point process
of excursions can be extended to include these downward jumps from strict
minima as follows. Let $\mathbf{x}$ denote the path $\mathbf{x}(t)=x$ for
all $t\ge 0$ and let 
\begin{equation}
\tilde{\mathcal{E}} =\mathcal{E}\cup \{\mathbf{x}:x<0\}.
\end{equation}
Define 
\begin{equation}
\tilde e_t= 
\begin{cases}
e_t, & \text{ if } e_t\in \mathcal{E} \\ 
\mathbf{x}, & \text{ if } \Delta \widehat{L}^{-1}_t=0 \text{ and } \Delta X_{%
\widehat{L}^{-1}_t}=x<0.%
\end{cases}%
\end{equation}
Then $\{(t,\tilde e_t):\tilde e_t\in \tilde{\mathcal{E}}\}$ is a Poisson
point process with characteristic measure $\tilde n$ given by 
\begin{equation}
\tilde n(A)=\widehat{n}(A\cap\mathcal{E})+d_{{\widehat{L}}^{-1}}\Pi_X^-(\{x:%
\mathbf{x}\in A\}).
\end{equation}

The properties of Poisson point processes used below can be found in
Proposition 0.2 of \cite{bert}. 
For $\delta \geq 0$ let 
\begin{equation}
A_{\delta }=\{\epsilon \in \mathcal{E}:h(\epsilon )>\delta \}
\end{equation}%
and $A_{\delta }^{c}=\tilde{\mathcal{E}}\setminus A_{\delta }$. Set 
\begin{equation*}
T_{\delta }=\inf \{t:e_{t}\in A_{\delta }\},
\end{equation*}%
and 
\begin{equation}
h^{(\delta )}=h(e_{T_{\delta}}),\ Z^{(\delta )}=\widehat{H}_{T_{\delta}-}%
\text{ and }D^{(\delta )}=|e_{T_{\delta}}(\zeta )|.  \label{hZD}
\end{equation}%
The case $\delta =0$ will only be considered when $|\widehat{n}|<\infty $.
Since the Poisson point processes $\{(t,e_{t}):e_{t}\in A_{\delta }\}$ and $%
\{(t,\tilde{e}_{t}):\tilde{e}_{t}\in A_{\delta }^{c}\}$ are independent, we
can write $\widehat{H}$ as the sum of two independent subordinators $%
\widehat{H}=J^{(\delta )}+K^{(\delta )}$ where 
\begin{equation}
K_{t}^{(\delta )}=\sum_{s\leq t}|e_{s}(\zeta )|I(e_{s}\in A)
\end{equation}%
is the sum of the jumps of $\widehat{H}$ that correspond to the ends of
excursions for which $h>\delta $, and $J^{(\delta )}=\widehat{H}-K^{(\delta
)} $. Their Laplace exponents are given by 
\begin{equation}
\begin{aligned} \kappa^{J^{(\gd)}} (\lambda) &= \widehat{d}
\lambda+\int_0^\infty (1-e^{-\gl x})\left\{\widehat{n}(h\le \delta,
|\epsilon(\zeta)|\in dx)+d_{{\whL}^{-1}}\Pi_X^-(dx)\right\},\\
\kappa^{K^{(\gd)}}(\lambda) &= \int_0^\infty (1-e^{-\gl x})\widehat{n}(h>
\delta, |\epsilon(\zeta)|\in dx) \end{aligned}
\end{equation}%
respectively. Here we are assuming $\widehat{q}=0$ which will be the case
below. Clearly $\widehat{H}_{t}=J_{t}^{(\delta )}$ for $t<T_{\delta }$ and $%
J^{(\delta )}$ does not jump at time $T_{\delta }$, so $Z^{(\delta
)}=J_{T_{\delta}}^{(\delta )}$. Further, $J^{(\delta )}$ is independent of $%
(T_{\delta},e_{T_{\delta}})$ and $e_{T_{\delta}}$ is independent of $%
T_{\delta}$, thus both $h^{(\delta )}$ and $D^{(\delta )}$ are independent
of $Z^{(\delta )}$. Additionally $T_{\delta}$ has an exponential
distribution with parameter $\widehat{n}(h>\delta )$, hence 
\begin{equation}
Ee^{-\lambda Z^{(\delta )}}=\int_{0}^{\infty }\widehat{n}(h>\delta )e^{-%
\widehat{n}(h>\delta )t}e^{-\kappa ^{J^{(\delta )}}(\lambda )t}dt=\frac{%
\widehat{n}(h>\delta )}{\widehat{n}(h>\delta )+\kappa ^{J^{(\delta
)}}(\lambda )}.  \label{ZL}
\end{equation}%
Since, by dominated convergence, 
\begin{equation}
\lim_{\delta \rightarrow 0}\kappa ^{J^{(\delta )}}(\lambda )=\widehat{d}%
\lambda +\int_{0}^{\infty }(1-e^{-\lambda x})\left\{ \widehat{n}%
(h=0,|\epsilon (\zeta )|\in dx)+d_{{\widehat{L}}^{-1}}\Pi
_{X}^{-}(dx)\right\} ,  \label{J0}
\end{equation}%
it follows from \eqref{ZL} that $Z^{(\delta )}\overset{\mathrm{P}}{%
\longrightarrow }0$ if either $\widehat{n}(h>0)=\infty $, or $\widehat{d}=0$%
, $\widehat{n}(h=0)=0$ and $d_{{\widehat{L}}^{-1}}=0$. Recall the condition $%
\widehat{n}(h=0)=0$ is equivalent to $X$ not being compound Poisson.
\bigskip

\noindent\textbf{Proof of Theorem \ref{T1}}\ \ Assume \eqref{tr}. We need to
consider three cases.

\noindent Case I: $\widehat{n}(h>0)=\infty$, or $\widehat{d}=d_{{\widehat{L}}%
^{-1}}=0$ and $X$ is not compound Poisson.

Recalling \eqref{hZD}, for any $x>\delta>0$ we have 
\begin{equation}
P(\tau_x<\infty)=P(h^{(\delta)}>x+Z^{(\delta)})+\int_0^\infty
P(h^{(\delta)}\le x+Z^{(\delta)}, Z^{(\delta)}+D^{(\delta)}\in
dy)P(\tau_{x+y}<\infty).
\end{equation}
Dividing by $P(\tau_x<\infty)$ and taking limits gives 
\begin{equation}
\begin{aligned} \lim_{x\to\infty}
\frac{P(h^{(\gd)}>x+Z^{(\gd)})}{P(\tx<\infty)}&=1-\int_0^\infty e^{-\ga
y}P(Z^{(\gd)}+D^{(\gd)}\in dy)\\ &=E(1-e^{-\ga (Z^{(\gd)}+D^{(\gd)})}).
\end{aligned}
\end{equation}
Since $h^{(\delta)}$ and $Z^{(\delta)}$ are independent and $h^{(\delta)}$
has distribution given by 
\begin{equation}
P(h^{(\delta)}\in \cdot\ )=\frac{\widehat{n}(h\in \cdot\ , h>\delta)}{%
\widehat{n}(h>\delta)},
\end{equation}
it then follows that 
\begin{equation}
\begin{aligned} \lim_{x\to\infty}
\frac{E[\whn(h>x+Z^{(\gd)})]}{\whn(h>\gd)P(\tx<\infty)}=E(1-e^{-\ga
(Z^{(\gd)}+D^{(\gd)})}). \end{aligned}
\end{equation}
Now for any $c>0$ 
\begin{equation}
\widehat{n}(h>x+c)P(Z^{(\delta)}\le c)\le E[\widehat{n}(h>x+Z^{(\delta)})]%
\le \widehat{n}(h>x),
\end{equation}
hence 
\begin{equation}
\begin{aligned} \widehat{n}(h>\delta) E(1-e^{-\ga(Z^{(\gd)}+D^{(\gd)})})&\le
\liminf_{x\to\infty}\frac{\whn(h>x)}{P(\tx<\infty)}\\ &\le
\limsup_{x\to\infty}\frac{\whn(h>x)}{P(\tx<\infty)} \le \frac{e^{\ga
c}\whn(h>\gd)}{P(Z^{(\gd)}\le c)}E(1-e^{-\ga(Z^{(\gd)}+D^{(\gd)})})
\end{aligned}
\end{equation}
where the last inequality again uses \eqref{tr}. Let $\delta\to 0$ then $%
c\to 0$ to obtain 
\begin{equation}
\begin{aligned} \limsup_{\gd\to
0}\widehat{n}(h>\delta)E(1-e^{-\ga(Z^{(\gd)}+D^{(\gd)})})&\le
\liminf_{x\to\infty}\frac{\whn(h>x)}{P(\tx<\infty)}\\ &\le
\limsup_{x\to\infty}\frac{\whn(h>x)}{P(\tx<\infty)}\\ &\le \liminf_{\gd\to
0}\widehat{n}(h>\delta)E(1-e^{-\ga(Z^{(\gd)}+D^{(\gd)})}). \end{aligned}
\end{equation}
Thus both limits exist and 
\begin{equation}  \label{limit}
\lim_{x\to\infty}\frac{\widehat{n}(h>x)}{P(\tau_x<\infty)}=\lim_{\delta\to 0}%
\widehat{n}(h>\delta)E(1-e^{-\alpha(Z^{(\delta)}+D^{(\delta)})}).
\end{equation}

To evaluate the limit observe that since $Z^{(\delta)}$ and $D^{(\delta)}$
are independent 
\begin{equation}
E(1-e^{-\alpha(Z^{(\delta)}+D^{(\delta)})})=E(1-e^{-\alpha
Z^{(\delta)}})+E(1-e^{-\alpha D^{(\delta)}}) -E(1-e^{-\alpha
Z^{(\delta)}})E(1-e^{-\alpha D^{(\delta)}}).
\end{equation}
By \eqref{ZL} and \eqref{J0} 
\begin{equation}
\begin{aligned} \lim_{\gd\to 0}\widehat{n}(h>\delta)E(1-e^{-\ga
Z^{(\gd)}})&=\lim_{\gd\to 0}\frac{ \whn(h> \gd)\gk^{J^{(\gd)}}(\ga)}{
\whn(h> \gd)+\gk^{J^{(\gd)}}(\ga)}\\ &=\widehat{d} \alpha+\int_0^\infty
(1-e^{-\ga z})d_{{\whL}^{-1}}\Pi_X^-(dz). \end{aligned}
\end{equation}
Next, since 
\begin{equation}  \label{Ddist}
P(D^{(\delta)}\in dz)=\frac{\widehat{n}(h> \delta, |\epsilon(\zeta)|\in dz)}{%
\widehat{n}(h>\delta)},
\end{equation}
we have by monotone convergence 
\begin{equation}
\begin{aligned}\label{D0} \widehat{n}(h> \delta)E(1-e^{-\ga
D^{(\gd)}})&=\widehat{n}(h> \delta)\left(1-\int_0^\infty e^{-\ga
z}\frac{\whn(h> \gd, |\gep(\zeta)|\in dz)}{\whn(h>\gd)}\right)\\
&=\int_0^\infty (1-e^{-\ga z})\widehat{n}(h> \delta, |\epsilon(\zeta)|\in
dz)\\ &\to \int_0^\infty (1-e^{-\ga z})\widehat{n}(|\epsilon(\zeta)|\in
dz).\\ \end{aligned}
\end{equation}
Finally by \eqref{D0} and $Z^{(\delta)}\overset{\mathrm{P}}{\longrightarrow}
0$, 
\begin{equation}
\widehat{n}(h>\delta)E(1-e^{-\alpha Z^{(\delta)}})E(1-e^{-\alpha
D^{(\delta)}})\to 0.
\end{equation}
Since \eqref{tr} implies $X_t\to -\infty$ a.s., this means $\widehat{q}=0$
and so by \eqref{piHn}, the limit in \eqref{limit} is $\widehat{\kappa}%
(\alpha)$. This proves \eqref{ter} which in turn implies \eqref{er}. \medskip

\noindent Case II: $\widehat{n}(h> 0)<\infty$ and $\widehat{d}> 0$ or $d_{{%
\widehat{L}}^{-1}}>0$.

%Note that the second condition implies that $X$ is not compound Poisson. 
If $0$ is irregular for $(-\infty ,0)$, then $(\widehat{L}^{-1},\widehat{H})$ is
bivariate compound Poisson, so $\widehat{d}=d_{{\widehat{L}}^{-1}}=0$. Thus
by Proposition \ref{P11}, it is necessarily the case that $0$ is irregular
for $[0,\infty )$ and $\overline{\Pi }_{X}^{+}(0)<\infty $. Hence $%
X_{t}=Y_{t}-U_{t}$ where $Y$ is a spectrally positive compound Poisson
process and $U$ is a subordinator which is not compound Poisson. The Laplace
exponent of $U$ is 
\begin{equation}
\kappa ^{U}(\lambda )=d_{U}\lambda +\int_{0}^{\infty }(1-e^{-\lambda x})\Pi_X
^{-}(dx).
\end{equation}%
Since $0$ is irregular for $[0,\infty )$, it suffices to prove the result
when $\widehat{L}$ is given by 
\begin{equation}
\widehat{L}_{t}=\int_{0}^{t}I(X_{s}=\setbox0=\hbox{$X$}\slantmathcorr=\wd0%
\underline{\hbox to 0.8\wd0{\vphantom{\hbox{$X$}}}}\hskip-0.8\wd0\hbox{$X$}%
_{s})ds.
\end{equation}%
In this case we have $\widehat{L}_{t}^{-1}=t$ until the time of the first
jump of $Y$, at which time $\widehat{L}^{-1}$ also jumps. Thus $d_{{\widehat{%
L}}^{-1}}=1$, $\widehat{d}=d_{U}$ and 
\begin{equation}
T=\inf \{t:\Delta Y_{t}>0\}=\inf \{t:e_{t}\in \mathcal{E}\}
\end{equation}%
has an exponential distribution with parameter $|\Pi ^{+}|=|\widehat{n}|.$
Setting $\delta =0$ in the discussion preceding the proof of Theorem \ref{T1}%
, we can write $\widehat{H}=J+K$ where $J$ has the same distribution as $U$
(since $d_{{\widehat{L}}^{-1}}=1$) and 
\begin{equation}
\begin{aligned} \kappa^K (\lambda) &= \int_0^\infty (1-e^{-\gl
x})\widehat{n}(|\epsilon(\zeta)|\in dx). \end{aligned}
\end{equation}%
Observe that $J\neq U$ but $J_{s}=U_{s}$ for $s\leq T$. Let $h_{1}=h(e_{T})$
be the height and $D_{1}=|e_{T}(\zeta )|$ the overshoot of the first
excursion. Then as in \eqref{D0} 
\begin{equation}
\begin{aligned}\label{D1} E(1-e^{-\gl D_1}) &=\frac{1}{|\whn|}\int_0^\infty
(1-e^{-\gl z})\widehat{n}(|\epsilon(\zeta)|\in dz)=\frac{\gk^K
(\gl)}{|\whn|} \end{aligned}
\end{equation}%
Again, as noted previously, $J$ is independent of $(T,e_{T})$ (this would
not be true if $J$ were replaced by $U$), and $e_{T}$ is independent of $T$.
In what follows, it will sometimes be convenient to write $P(\tau
_{x}<\infty )$ as $P(\overline{X}_{\infty }>x)$ where $\overline{X}%
_{t}=\sup_{s\leq t}X_{s}.$ We also write $S$ for the righthand endpoint $%
\widehat{L}_{T}^{-1}$ of the first excursion interval.  
Then for any $t>0$, since $J_{s}=U_{s}$ for $s\leq T$, 
\begin{equation}
\begin{aligned} P(\tau_x<\infty)&=P(T\le t, h_1> x+J_T)+P(T\le t, h_1\le
x+J_T, \sup_{r\ge 0}(X_{S+r}-X_S)>x+J_T+D_1)\\ &\qquad\qquad +P(T>t,
\sup_{r\ge 0}(X_{t+r}-X_t)>x+J_t)\\ &=\int_0^t P(T\in
ds)P(h_1>x+J_s)+\int_0^t P(T\in ds)Ef_x(h_1,D_1,J_s)+P(T>t)Eg_x(J_t)\\
\end{aligned}
\end{equation}%
where 
\begin{equation}
f_{x}(y,z,w)=I(y\leq x+w)P(\overline{X}_{\infty }>x+w+z),\ \ g_{x}(w)=P(%
\overline{X}_{\infty }>x+w)
\end{equation}%
Thus dividing by $P(\tau _{x}<\infty )$ and letting $x\rightarrow \infty $,
we obtain 
\begin{equation}
\begin{aligned} \lim_{x\to\infty}\int_0^t P(T\in
ds)\frac{P(h_1>x+J_s)}{P(\tx<\infty)}&=1-P(T>t)Ee^{-\ga J_t}-\int_0^t P(T\in
ds)Ee^{-\ga(J_s+D_1)}\\ &=1-e^{-(|\whn|+\gk^J(\ga))t} -
|\widehat{n}|Ee^{-\ga D_1}\int_0^t e^{-(|\whn|+\gk^J(\ga))s}ds. \end{aligned}
\end{equation}%
Now divide by $t$ and let $t\rightarrow 0$ to get 
\begin{equation}
\begin{aligned} \lim_{t\to 0}\lim_{x\to\infty}\int_0^t P(T\in
ds)\frac{P(h_1>x+J_s)}{tP(\tx<\infty)} &=|\widehat{n}|+\kappa^J(\alpha) -
|\widehat{n}|\left(1-\frac{\gk^K(\ga)}{|\whn|}\right)=\widehat{\kappa}(%
\alpha). \end{aligned}
\end{equation}%
Since 
\begin{equation}
P(h_{1}>x+J_{t})\leq P(h_{1}>x+J_{s})\leq P(h_{1}>x)
\end{equation}%
for $0\leq s\leq t$ and $J_{t}\overset{\mathrm{P}}{\longrightarrow }0$, it
then easily follows that 
\begin{equation}
\lim_{x\rightarrow \infty }\frac{|\widehat{n}|P(h_{1}>x)}{P(\tau _{x}<\infty
)}=\widehat{\kappa }(\alpha )
\end{equation}%
which is equivalent to \eqref{ter}.

\medskip

\noindent Case III: $X$ is compound Poisson.

This reduces to the random walk case. Dividing through \eqref{b} by $%
P(\tau_x<\infty)$ and letting $x\to \infty$ gives 
\begin{equation}  \label{CPcase}
\lim_{x\to\infty}\frac{P(h_1>x)}{P(\tau_x<\infty)}=E(1-e^{-\alpha |S_{%
\widehat{T}_1}|}).
\end{equation}
Since $\widehat{H}_1$ is the sum of an independent Poisson, with parameter $|%
\widehat{n}|$, number of IID copies of $|S_{\widehat{T}_1}|$, it easily
follows that the limit in \eqref{CPcase} is $\widehat{\kappa}(\alpha)/|%
\widehat{n}|$ which is equivalent to \eqref{ter}. \medskip

In the converse direction, assume \eqref{er}. 
By the compensation formula 
\begin{equation}
\begin{aligned}\label{t/n}
P(\tau_x<\infty)&=E\sum_{t}I(\overline{X}_{\whL_{t-}^{-1}}\le x,
h(e_t)>x+|\setbox0=\hbox{$X$} \slantmathcorr=\wd0 \underline{\hbox to
0.8\wd0{\vphantom{\hbox{$X$}}}} \hskip-0.8\wd0\hbox{$X$}
_{\whL_{t-}^{-1}}|)\\ &=E\int_0^\infty dt I(\overline{X}_{\whL_{t-}^{-1}}\le
x) \widehat{n}(h>x+|\setbox0=\hbox{$X$} \slantmathcorr=\wd0 \underline{\hbox
to 0.8\wd0{\vphantom{\hbox{$X$}}}} \hskip-0.8\wd0\hbox{$X$}
_{\whL_{t-}^{-1}}|)\\ &=\int_0^\infty dt \int_{y\ge
0}P(\overline{X}_{\whL_{t-}^{-1}}\le x, |\setbox0=\hbox{$X$}
\slantmathcorr=\wd0 \underline{\hbox to 0.8\wd0{\vphantom{\hbox{$X$}}}}
\hskip-0.8\wd0\hbox{$X$} _{\whL_{t-}^{-1}}|\in dy)\widehat{n}(h>x+y).
\end{aligned}
\end{equation}
By \eqref{Pot}, for any $\varepsilon\in (0,\alpha)$ there exists a constant $%
A$ such that 
\begin{equation}
\frac{\widehat{n}(h>x+y)}{\widehat{n}(h>x)}\le A e^{-(\alpha-\varepsilon)
y}\ \text{ for all } x\ge 1, y\ge 0.
\end{equation}
Thus for $x\ge 1$ 
\begin{equation}
\begin{aligned} \int_0^\infty dt \int_{y\ge
0}P(\overline{X}_{\whL_{t-}^{-1}}\le x, &|\setbox0=\hbox{$X$}
\slantmathcorr=\wd0 \underline{\hbox to 0.8\wd0{\vphantom{\hbox{$X$}}}}
\hskip-0.8\wd0\hbox{$X$} _{\whL_{t-}^{-1}}|\in
dy)\frac{\whn(h>x+y)}{\whn(h>x)}\\ &\le A\int_0^\infty dt \int_{y\ge 0}P(
|\setbox0=\hbox{$X$} \slantmathcorr=\wd0 \underline{\hbox to
0.8\wd0{\vphantom{\hbox{$X$}}}} \hskip-0.8\wd0\hbox{$X$}
_{\whL_{t-}^{-1}}|\in dy)e^{-(\ga-\vep) y}\\ &\le A\int_0^\infty dt E
e^{-(\ga-\vep) \whH_{t}}\\ &=\frac{A}{\whk(\ga-\vep)}<\infty. \end{aligned}
\end{equation}
Hence, dividing \eqref{t/n} by $\widehat{n}(h>x)$ and applying dominated
convergence we obtain 
\begin{equation}
\begin{aligned} \lim_{x\to\infty} \frac{P(\tx<\infty)}{\whn(h>x)}
&=\int_0^\infty dt \int_{y\ge 0}P( |\setbox0=\hbox{$X$} \slantmathcorr=\wd0
\underline{\hbox to 0.8\wd0{\vphantom{\hbox{$X$}}}} \hskip-0.8\wd0\hbox{$X$}
_{\whL_{t-}^{-1}}|\in dy)e^{- \ga y} =\frac{1}{\whk(\ga)}. \end{aligned}
\end{equation}
Thus \eqref{ter} holds which in turn implies \eqref{tr}. \qquad$%
\sqcup\!\!\!\!\sqcap$ \bigskip

\begin{remark}
\label{Rem0} If $\Gamma^*=0$ in \eqref{h}, then a simpler version of the
above proof where dividing by $P(\tau_x<\infty)$ is replaced by dividing by $%
e^{-\gamma x}$, shows that the limit in \eqref{g} is also $0$.
\end{remark}

\begin{remark}
\label{Rem1} If $X_t\not\to -\infty$ a.s. then $P(\tau_x<\infty)= 1$ for all 
$x\ge 0$, so \eqref{tr} trivially holds when $\alpha=0$. Since this provides
no useful information about the asymptotic behaviour of $P(\tau_x<\infty)$,
we must also include the condition $X_t\to -\infty$ when considering %
\eqref{tr} in the $\alpha=0$ case. In that case the proof for $\alpha>0$ is
easily modified, and is in fact much simpler, to show that \eqref{ter} holds
with $\alpha =0$, the limit being $\widehat{\kappa}(0)=0$ since $X_t\to
-\infty$. However this does not enable us to conclude anything about %
\eqref{er}. Conversely if \eqref{er} holds with $\alpha=0$, then we can
divide through \eqref{t/n} by $\widehat{n}(h>x)$ and apply Fatou to obtain 
\begin{equation}
\begin{aligned} \liminf_{x\to\infty} \frac{P(\tx<\infty)}{\whn(h>x)}
&\ge\int_0^\infty dt \int_{y\ge 0}P( |\setbox0=\hbox{$X$}
\slantmathcorr=\wd0 \underline{\hbox to 0.8\wd0{\vphantom{\hbox{$X$}}}}
\hskip-0.8\wd0\hbox{$X$} _{\whL_t^{-1}}|\in dy)=\widehat{V}(\infty)=\frac
1{\whk(0)}. \end{aligned}
\end{equation}
The corresponding upper bound holds trivially for every $x>0$ without taking
the limit. Thus \eqref{ter} holds with $\alpha =0$, but we are unable to
conclude anything about \eqref{tr} unless $\widehat{q}>0$. In this direction
there is no need to assume $X_t\to -\infty$ a.s. If $X_t\not\to -\infty$
a.s. then \eqref{ter} simply reduces to $\widehat{n}(h=\infty)=\widehat{q}$.
\end{remark}
\bigskip

\noindent\textbf{Proof of Theorem \ref{KKMC}.}\ \ Assume $\Pi_{H}\in 
\mathcal{L}^{(\alpha)}$. By Vigon's \'equation amicale, see (5.3.3) of \cite%
{D}, for any $t>0$, 
\begin{equation}
\overline{\Pi}_X(t)=\int_0^\infty \Pi_{H}(t+dy)\overline{\Pi}_{\widehat{H}%
}(y)+\widehat{d} \Pi_{H}^{\prime}(t)+\widehat{q} \overline{\Pi}_{H}(t),
\end{equation}
where $\Pi_{H}^{\prime}$ denotes the cadlag version of the density of $\Pi_H$%
, which exists when $\widehat{d}>0$. By Fubini's Theorem 
\begin{equation}
\overline{\Pi}_X(t)=\int_0^\infty (\overline{\Pi}_{H}(t)-\overline{\Pi}%
_{H}(t+y))\Pi_{\widehat{H}}(dy)+\widehat{d} \Pi_{H}^{\prime}(t)+\widehat{q} 
\overline{\Pi}_{H}(t),
\end{equation}
thus 
\begin{equation}
\begin{aligned}\label{XH1} \frac 1{\pibar_{H}(x)}\int_x^\infty
\overline{\Pi}_X(t)dt =\int_0^\infty
\Pi_{\whH}(dy)\int_0^y\frac{\pibar_{H}(x+t)}{\pibar_{H}(x)}dt +\widehat{d}+
\frac {\whq}{\pibar_{H}(x)}\int_x^\infty \overline{\Pi}_{H}(t)dt.
\end{aligned}
\end{equation}
Fix $\varepsilon\in (0, \alpha)$. By \eqref{Pot}, for some $A$ and all $x\ge
1, y\ge 0$ 
\begin{equation}  \label{XH0}
\int_0^y\frac{\overline{\Pi}_{H}(x+t)}{\overline{\Pi}_{H}(x)}dt\le A\int_0^y
e^{-(\alpha-\varepsilon)t} dt =\frac{A(1-e^{-(\alpha-\varepsilon)y})}{%
\alpha-\varepsilon}.
\end{equation}
This final expression is integrable over $(0,\infty)$ with respect to $\Pi_{%
\widehat{H}}(dy)$, hence we may apply dominated convergence to conclude 
\begin{equation}  \label{XH2}
\int_0^\infty \Pi_{\widehat{H}}(dy)\int_0^y\frac{\overline{\Pi}_{H}(x+t)}{%
\overline{\Pi}_{H}(x)}dt \to \int_0^\infty \Pi_{\widehat{H}}(dy) \frac{%
(1-e^{- \alpha y})}{\alpha}.
\end{equation}
Similarly, another appeal to \eqref{Pot} together with dominated convergence gives
\begin{equation}
\begin{aligned}\label{XH3} \frac {\whq}{\pibar_{H}(x)}\int_x^\infty
\overline{\Pi}_{H}(t)dt&=\widehat{q}\int_0^\infty\frac{\pibar_{H}(x+t)}{%
\pibar_{H}(x)}dt \to \frac{\whq}{\ga}. \end{aligned}
\end{equation}
Thus by \eqref{XH1}, \eqref{XH2} and \eqref{XH3} 
\begin{equation}
\frac 1{\overline{\Pi}_{H}(x)}\int_x^\infty \overline{\Pi}_X(t)dt\to\frac{%
\widehat{\kappa}(\alpha)}{\alpha}.
\end{equation}
Now fix $a>0$. Then 
\begin{equation}
\begin{aligned} \frac{a\pibar_X(x)}{\pibar_{H}(x)}&\le \frac 1
{\pibar_{H}(x)}\int_{x-a}^x \overline{\Pi}_X(t)dt\\ &=
\frac{\pibar_{H}(x-a)}{\pibar_{H}(x)}\frac 1
{\pibar_{H}(x-a)}\int_{x-a}^\infty \overline{\Pi}_X(t)dt - \frac 1
{\pibar_{H}(x)}\int_{x}^\infty \overline{\Pi}_X(t)dt\\ &\to
\frac{\whk(\ga)}{\ga}(e^{\ga a}-1). \end{aligned}
\end{equation}
Divide by $a$ and let $a\to 0$ to obtain 
\begin{equation}
\limsup_{x\to\infty} \frac{\overline{\Pi}_X(x)}{\overline{\Pi}_{H}(x)}\le 
\widehat{\kappa}(\alpha).
\end{equation}
Integrating over $[x,x+a]$ gives the corresponding lower bound. Hence %
\eqref{piXHr} holds and consequently $\overline{\Pi}_{X}\in \mathcal{L}%
^{(\alpha)}$.

The opposite direction is straightforward. Assume $\overline{\Pi}_{X}\in 
\mathcal{L}^{(\alpha)}$. By Vigon's equation amicale invers\'ee, see (5.3.4)
of \cite{D}, for $x>0$ 
\begin{equation}  \label{VA}
\frac{\overline{\Pi}_{H}(x)}{\overline{\Pi}_X(x)}=\int_0^\infty \widehat{V}%
(dy)\frac{\overline{\Pi}_X(x+y)}{\overline{\Pi}_X(x)}
\end{equation}
To take the limit inside the integral, we again we use \eqref{Pot} and
observe 
\begin{equation}
\int_0^\infty \widehat{V}(dy)Ae^{-(\alpha-\varepsilon)y} = \int_0^\infty 
\widehat{V}(y)A(\alpha-\varepsilon) e^{-(\alpha-\varepsilon)y} dy<\infty
\end{equation}
since $\widehat{V}(y)\le Cy$ for $y\ge 1$ by Proposition III.1 of \cite{bert}%
. Thus by dominated convergence 
\begin{equation}
\frac{\overline{\Pi}_{H}(x)}{\overline{\Pi}_X(x)}\to \int_0^\infty \widehat{V%
}(dy)e^{-\alpha y}=\frac{1}{\widehat{\kappa}(\alpha)}.
\end{equation}
Hence \eqref{piXHr} holds and consequently $\overline{\Pi}_{H}\in \mathcal{L}%
^{(\alpha)}$. \qquad$\sqcup\!\!\!\!\sqcap$ \bigskip

\begin{remark}
When $\alpha=0$, \eqref{piXr} implies \eqref{piXHr} and \eqref{piHr} implies %
\eqref{piXHr}, but \eqref{piXr} and \eqref{piXHr} are not necessarily
equivalent since it is possible that $\widehat{\kappa}(0)=0$. To see this,
by \eqref{VA} for any $x>0$ without any assumptions on $\overline{\Pi}_X$ or 
$\overline{\Pi}_H$, 
\begin{equation}
\frac{\overline{\Pi}_{H}(x)}{\overline{\Pi}_X(x)}\le \widehat{V}%
(\infty)=\frac 1{\widehat{\kappa}(0)}.
\end{equation}
If $\overline{\Pi}_X\in \mathcal{L}^{(\alpha)}$, then applying Fatou to %
\eqref{VA} proves \eqref{piXHr}. If $\overline{\Pi}_H\in \mathcal{L}%
^{(\alpha)}$ then for any $K>0$ 
\begin{equation}
\frac{\overline{\Pi}_{H}(x)}{\overline{\Pi}_H(x+K)}\ge\int_0^K \widehat{V}%
(dy)\frac{\overline{\Pi}_X(x+y)}{\overline{\Pi}_H(x+K)} \ge \widehat{V}(K)%
\frac{\overline{\Pi}_X(x+K)}{\overline{\Pi}_H(x+K)}.
\end{equation}
Letting $x\to\infty$ and then $K\to\infty$ proves \eqref{piXHr}.
\end{remark}

\noindent\textbf{Proof of Theorem \ref{LPP}}\ \ We prove the equivalent
formulation in Remark \ref{Rem4}. Thus assume $\overline{\Pi }_{H}\in 
\mathcal{L}^{(\alpha)}$ , $E(e^{\alpha H_{1}})<1,$ and \eqref{tHrat} hold.
Let $Z= H_{\tau_1}$ if $\tau_1<\infty$ and set $Z$ equal to some cemetery
state otherwise. Then by Proposition III.2 of \cite{bert}, 
\begin{equation}  \label{HLa1}
\lim_{x\to\infty}\frac{P(Z>x)}{\overline{\Pi}_{H}(x)}=\lim_{x\to\infty}%
\int_0^1 V(dz)\frac{\overline{\Pi}_{H}(x-z)}{\overline{\Pi}_{H}(x)}=
\int_0^1 e^{\alpha z}V(dz).
\end{equation}
Hence $P(Z>x)\in \mathcal{L}^{(\alpha)}$. Further, since $Ee^{\alpha H_1}<
\infty$ implies $\int_1^\infty e^{\alpha y}\Pi_H (dy)< \infty$ by Theorem
25.17 of \cite{sato}, which in turn is equivalent to 
\begin{equation}  \label{PiHint}
\int_1^\infty \overline{\Pi}_H (y) \alpha e^{\alpha y}dy< \infty,
\end{equation}
we have 
\begin{equation}  \label{EZf}
Ee^{\alpha Z}=\int_0^\infty P(Z>y) \alpha e^{\alpha y}dy< \infty.
\end{equation}
Now for $x>1$, 
\begin{equation}
\begin{aligned}\label{E0}
\frac{P(\tx^H<\infty)}{\pibar_{H}(x)}=\frac{P(Z>x)}{\pibar_{H}(x)}+\int_0^x
P(Z\in dy)\frac{P(\tau^H_{x-y}<\infty)}{\pibar_{H}(x)}. \end{aligned}
\end{equation}
By bounded convergence, then monotone convergence 
\begin{equation}
\begin{aligned}\label{0Klim} \lim_{K\to\infty}\lim_{x\to\infty}\int_{0}^K
P(Z\in dy)\frac{P(\tau^H_{x-y}<\infty)}{\pibar_{H}(x)}
&=\lim_{K\to\infty}\int_0^K P(Z\in dy) L'e^{\ga y}\\ &=\int_0^\infty P(Z\in
dy) L'e^{\ga y}<\infty \end{aligned}
\end{equation}
by \eqref{EZf}, while 
\begin{equation}
\begin{aligned}\label{E00} \lim_{K\to\infty}\lim_{x\to\infty}\int_{x-K}^x
&P(Z\in dy)\frac{P(\tau^H_{x-y}<\infty)}{\pibar_{H}(x)}\\ &=
\lim_{K\to\infty}\lim_{x\to\infty}\int_0^K P(\tau^H_{y}<\infty)\frac{P(Z\in
x-dy)}{P(Z>x)}\frac{P(Z>x)}{\pibar_{H}(x)}\\ &= \lim_{K\to\infty}\int_0^K
P(\tau^H_{y}<\infty)\alpha e^{\ga y} dy \int_0^1 e^{\ga z}V(dz)\\
&=\int_0^\infty P(\tau^H_{y}<\infty)\alpha e^{\ga y} dy \int_0^1 e^{\ga
z}V(dz)<\infty \end{aligned}
\end{equation}
by \eqref{tHrat} and \eqref{PiHint}. Thus by \eqref{tHrat}, \eqref{E0}, %
\eqref{0Klim} and \eqref{E00} 
\begin{equation}
\lim_{K\to\infty}\lim_{x\to\infty}\int_K^{x-K} P(Z\in dy)\frac{%
P(\tau^H_{x-y}<\infty)}{\overline{\Pi}_{H}(x)}\ \text{ exists.}
\end{equation}
By \eqref{tHrat} and \eqref{HLa1} it then follows that 
\begin{equation}  \label{E2}
\lim_{K\to\infty}\lim_{x\to\infty}\int_K^{x-K} P(Z\in dy)\frac{P(Z>x-y)}{%
P(Z>x)}\ \text{ exists.}
\end{equation}
Now let $Z_1,Z_2$ be independent and distributed as $Z$ conditional on $%
\tau_1<\infty$, so $P(Z_i\in dy)= P(Z\in dy)/P(\tau_1<\infty)$ for $i=1,2$.
Then $Z_1\in \mathcal{L}^{(\alpha)}$ and by (1) of \cite{EG} 
\begin{equation}
\begin{aligned} \int_K^{x-K}P(Z_1\in dy)\frac{P(Z_1> x-y)}{P(Z_1>x)}
&=\frac{P(Z_1+Z_2>x)}{P(Z_1>x)} -2\int_0^K\frac{P(Z_1>
x-y)}{P(Z_1>x)}P(Z_1\in dy)\\ &\qquad-\frac{P(Z_1>x-K)P(Z_1>K)}{P(Z_1>x)}\\
&\sim \frac{P(Z_1+Z_2>x)}{P(Z_1>x)}-2\int_0^K e^{\ga y}P(Z_1\in dy)-e^{\ga
K}P(Z_1>K)\\ \end{aligned}
\end{equation}
as $x\to\infty$. Thus letting $K\to \infty$ 
we obtain 
\begin{equation}
\begin{aligned} \lim_{x\to\infty}\frac{P(Z_1+Z_2>x)}{P(Z_1>x)}=
\lim_{K\to\infty}\lim_{x\to\infty}\int_K^{x-K} P(Z_1\in
dy)\frac{P(Z_1>x-y)}{P(Z_1>x)}+2Ee^{\ga Z_1} \end{aligned}
\end{equation}
which exists by \eqref{EZf} and \eqref{E2}. 
This means that $Z_1\in \mathcal{S}^{(\alpha)}$ and since $\mathcal{S}%
^{(\alpha)}$ is closed under tail equivalence, this in turn implies $%
\overline{\Pi}_{H}\in \mathcal{S}^{(\alpha)}$.

The converse holds by Lemma 3.5 of \cite{KKM}, where the value of $L^{\prime}$
is also calculated. \qquad$\sqcup\!\!\!\!\sqcap$

\begin{remark}
The equivalent formulation of Theorem \ref{LPP}, given in Remark \ref{Rem4},
continues to hold when $\alpha=0$ with the interpretation that $Ee^{\alpha
H_1}<1$ means $H$ is defective, i.e. $X_t\to -\infty$. The proof is an
obvious modification of the proof in the $\alpha>0$ case. Theorem \ref{LPP}
as stated does not hold for $\alpha=0$. This is because when $\overline{\Pi}%
_X\in\mathcal{L}^{(0)}$ one can show 
\begin{equation}
\liminf_{x\to\infty}\frac{\widehat{n}(h>x)}{\overline{\Pi}_X(x)}\ge \frac 1q.
\end{equation}
Thus if in addition \eqref{limXa0} holds, then $P(\tau_x<\infty)\in \mathcal{%
L}^{(0)}$ and $X_t\to-\infty$, hence by Remark \ref{Rem1} 
\begin{equation}
\lim_{x\to\infty}\frac{P(\tau_x<\infty)}{\overline{\Pi}_X(x)}=\infty.
\end{equation}
\end{remark}

\end{document}